\documentstyle[12pt]{article}

\begin{document}
\hbadness=10000
\hbadness=10000
\begin{titlepage}
\nopagebreak
%\vspace*{1.5 cm}
\def\thefootnote{\fnsymbol{footnote}}
%%%%%%%%%%%%%%%%%%%%%%%%%% Preprint No. %%%%%%%%%%%
%\rightline{}
\begin{flushright}
        {\normalsize
 LMU-TPW 94-27 \\
q-alg 9501028 \\
Mod. Phys. Lett. A }\\
\end{flushright}
%%%%%%%%%%%%%%%%%%%%%%%%%%%%%%%%%%%%%%%%%%%%%%%%%%%
\vspace{2cm}
\vfill
\begin{center}
\renewcommand{\thefootnote}{\fnsymbol{footnote}}
{\large \bf Quantum Deformation of $igl(n)$ Algebra on Quantum Space}

\vfill
\vspace{1.1cm}

{\bf Tatsuo Kobayashi 
\footnote[1]{Alexander von Humboldt Fellow \\
\phantom{xxx}e-mail:kobayash@lswes8.ls-wess.physik.uni-muenchen.de}
}

\vspace{0.5cm}
       Sektion Physik, Universit\"at M\"unchen, \\

       Theresienstrasse 37, D-80333 M\"unchen, Germany \\
              
and 

{\bf Haru-Tada Sato 
\footnote[2]{Fellow of Danish Research Academy\\
\phantom{xxx} e-mail:sato@nbi.dk }}

\vspace{0.5cm}
       The Niels Bohr Institute, University of Copenhagen, \\

       Blegdamsvej 17, DK-2100 Copenhagen, Denmark \\

\vfill
\end{center}
\vspace{1.1cm}

\vfill
\nopagebreak
\begin{abstract}
We study quantum deformed $gl(n)$ and $igl(n)$ algebras on a quantum 
space discussing multi-parametric extension. 
We realize elements of deformed $gl(n)$ and $igl(n)$ algebras by 
a quantum fermionic space. 
We investigate a map between deformed $igl(2)$ algebras of our basis 
and other basis.

\end{abstract}

\vfill
\end{titlepage}
\pagestyle{plain}
\newpage
\voffset = 0.5 cm

%%%%%%%%%%%%%%%%%%%%%%%%%%% Section 1 %%%%%%%%%%%%%%%%%%%%%%%%%%%%%%%%%%%%%%
\vspace{0.8 cm}
\leftline{\large \bf 1. Introduction}
\vspace{0.8 cm}
%%%%%%%%%%%%%%%%%%%%%%%%%%%%%%%%%%%%%%%%%%%%%%%%%%%%%%%%%%%%%%%%%%%%%%%%%%%%

Quantum deformation of groups and algebras has been studied in recent years 
[1-4]. In addition to deformation of Lie groups and algebras, some types of 
inhomogeneous groups and algebras are deformed, e.g., deformations of 
$igl(n)$, $iso(n)$ and Poincar\'e algebras [5-15].

Quantum spaces are non-commutative spaces to represent quantum groups and 
their commutation relations between coordinates and derivatives are 
governed by the R-matrices for corresponding quantum groups [16-18].
Further, one of quantum deformed $su(2)$ algebras is constructed based on 
the two-dimensional quantum space in ref.\cite{SWZ2}.
In ref.\cite{KU2} supersymmetric extension of the quantum space, i.e., 
quantum bosonic and fermionic spaces, is discussed.
The quantum $su(2)$ is extended into construction of quantum deformed 
$su(m|n)$ algebra on the quantum superspace in ref.\cite{Kobayashi}, which 
also shows that generators of the deformed $su(2)$ algebra can be realized in 
terms of the two-dimensional quantum fermionic coordinates and derivatives. 
In this paper, we hence generalize these construction of quantum algebras 
into an inhomogeneous case. The purpose of this paper is to explain how to 
construct quantum deformed algebra $igl(n)$ on fermionic quantum space as 
well as on bosonic one.

We employ the quantum Grassmann representation to construct quantum 
algebras. This method is used in some publications mentioned above. 
Grassmannian construction of quantum groups has a powerful advantage because 
its representation becomes very simple compared to usual construction of 
quantum groups. However, it is not clear whether or not the deformed 
algebras based on this method possess non-trivial Hopf algebra 
structure. None of the non-trivial structure has been found so far. 
It is hence important to find a connection of our approach to known Hopf 
algebra structures. In this sense, we should discuss also a relation 
between our deformed algebra and other deformed algebra which exhibits 
non-trivial Hopf algebra structure.                             

This paper is organized as follows. 
In section 2, we explain how to construct the deformed 
$su(n)$ algebras on the quantum spaces in preliminary. As new results,  
we mention a construction of the $su(2)$ algebra on quantum space 
with two deformation parameters and show a realization of deformed $su(n)$ 
in terms of fermionic coordinates and derivatives.  
In section 3, we define a dilatation generator $D$ in order to 
construct deformed $gl(n)$ and $igl(n)$ algebras.
In section 4, we exhibit a relation of our deformed $igl(2)$ algebra to 
Castellani's one. Section 5 is devoted to conclusion and discussion.

%%%%%%%%%%%%%%%%%%%%%%%%%%%% Section 2 %%%%%%%%%%%%%%%%%%%%%%%%%%%%%%%%%%%%%%
\vspace{0.8 cm}
\leftline{\large \bf 2. Deformed $su(n)$ algebra on quantum space}
\vspace{0.8 cm}
%%%%%%%%%%%%%%%%%%%%%%%%%%%%%%%%%%%%%%%%%%%%%%%%%%%%%%%%%%%%%%%%%%%%%%%%%%%%5

Quantum space is the non-commutative space representing a quantum group. 
Differential calculi on bosonic and fermionic quantum spaces for the 
quantum deformed $SU(n)$ are discussed in refs.\cite{WZ,KU2}. 
The quantum group $GL_q(n)$ has more deformation parameters than $SU_q(n)$ 
does. For example $GL_q(2)$ has two deformation parameters. 
Quantum space with multi-parameters were studied in ref.
\cite{SWZ3,Schirrmacher,KU2}.

First of all, we show a deformed $su(2)$ algebra on 2-dimensional quantum 
fermionic space with two deformation parameters as a simple example of 
multi-parameter cases. Bosonic case is remarked at the end of this section. 
We set up commutation relations between $\theta^\alpha$ and 
$\partial_\alpha$ with two parameters, $r$ and $q$, as follows,
$$
q\theta^1 \theta^2=-\theta^2 \theta^1, \quad 
\partial_1 \partial_2=-{r^2 \over q} \partial_2 \partial_1,$$
$$ q\partial_2 \theta^1 = - \theta^1 \partial_2, 
\quad , \partial_1 \theta^2 = - {q \over r^2}\theta^2 \partial_1, 
\eqno(2.1)$$
$$ \partial_\alpha \theta_\alpha = 1- \theta_\alpha \partial_\alpha 
+(r^{-2}-1)\sum_{\beta > \alpha} \theta ^\beta \partial _\beta,
$$
where $(\theta^\alpha)^2=(\partial_\alpha)^2=0$.
When we take $q=r$, we obtain the quantum plane with one parameter.

Commutation relations among generators are found as follows. First we 
consistently define operations of generators on the above quantum space. 
We then construct commutation relations among generators with the aid of 
the operations on the quantum space. Now let us consider two generators 
denoted by $T^1_2$ and $T^2_1$, which correspond to 
$\theta^1\partial_2$ and $\theta^2\partial_1$ in the classical limit 
($r,q \rightarrow 1$).
Following ref.\cite{SWZ2}, we assume that $T^1_2$ acts on $\theta^\alpha$ 
and $\partial_\alpha$ as follows,
$$
T^1_2 \theta^\alpha = a_\alpha \theta^\alpha T^1_2 +\delta^\alpha_2 
\theta^2,$$
$$T^1_2 \partial_\alpha=b_\alpha \partial_\alpha + 
c\delta^1_\alpha \partial_2.
\eqno(2.2)
$$
We investigate consistency between (2.1) and (2.2).
For example, we calculate $T^1_2 (\theta^2)^2$ as 
$$
T^1_2(\theta^2)^2=(a_2)^2(\theta^2)^2T^1_2+a_2\theta^2\theta^1+
\theta^1\theta^2.
\eqno(2.3)
$$
Compared with (2.1), we get $a_2=1/q$. Other relations are obtained in 
the same way
$$
[T^1_2,\theta^1]_{r^2/q}=0, \quad 
[T^1_2,\theta^2]_{1/q}=\theta^1,$$
$$[T^1_2,\partial_1]_{q/r^2}=-{q \over r^2}\partial_2, \quad 
[T^1_2,\partial_2]_{q}=0,
\eqno(2.4)
$$
where $[A,B]_h \equiv AB-hBA$.
Similarly for $T^2_1$, its operations on the quantum space are
$$
[T^2_1,\theta^1]_{q}=\theta^2, \quad 
[T^2_1,\theta^2]_{q/r^2}=0,$$
$$[T^2_1,\partial_1]_{1/q}=0, \quad 
[T^2_1,\partial_2]_{r^2/q}=-{r^2 \over q}\partial_1.
\eqno(2.5)
$$
As for a Cartan element $H_1$, it is defined by a commutation relation 
between $T^1_2$ and $T^2_1$. Since we have the following relations on the 
coordinates
$$
[T^2_1T^1_2,\theta^1]_{r^2}={r^2 \over q}\theta^2 T^1_2, \quad 
[T^1_2T^2_1,\theta^1]_{r^2}=q^{-1}\theta^2T^1_2 + \theta^1,
\eqno(2.6)
$$
it seems natural to define $H_1$ as 
$$
H_1 \equiv r^{-1}[T^2_1,T^1_2]_{r^2},
\eqno(2.7)
$$
in order to eliminate the linear terms of $T^1_2$ on the right hand side 
of (2.6). Accordingly, operations of $H_1$ on $\theta^\alpha$ and 
$\partial_\alpha$ 
become 
$$
[H_1 ,\theta^1]_{r^2}=-r\theta^1, \quad 
[H_1 ,\partial_1]_{1/r^2}=r^{-1}\partial_1 ,$$
$$[H_1 ,\theta^2]_{1/r^2}=r^{-1}\theta^2, \quad 
[H_1 ,\partial_2]_{r^2}=-r\partial_2.
\eqno(2.8)
$$

Using (2.4), (2.5) and (2.8), we have the following relations,
$$
[H_1,T^1_2]_{r^4}=-r^2(r+r^{-1})T^1_2,\quad
[H_1,T^2_1]_{1/r^4}=r^{-2}(r+r^{-1})T^2_1 ,
\eqno(2.9)
$$
Eqs.(2.7) and (2.9) compose a deformed $su(2)$ algebra as an abstract 
algebra. Note that this deformed $su(2)$ algebra on the quantum plane with 
the two parameters include only one parameter $r$ and the operations of 
$H_1$ is independent on $q$. The operations of only $T^\alpha_\beta$, (2.4) 
and (2.5), include the two parameters. For one parameter deformation, it was 
shown in ref.\cite{Kobayashi} that the deformed $su(2)$ generators are 
realized as
$$
T^1_2=\theta^1 \partial_2, \quad T^2_1=\theta^2 \partial_1 $$
$$H_1=r^{-1}\theta^2 \partial_2-r\theta^1 \partial_1.
\eqno(2.10)
$$
This is of course true even in the case with two deformation parameters.

In a way similar to the above derivation, we can construct a deformed 
$su(n)$ algebra using $n$-dimensional quantum space. We here show the 
results of a case with only one deformation parameter $r$, for simplicity. 
Multi-parametric extensions can be derived in a straightforward way. 
Commutation relations between fermionic coordinates $\theta^\alpha$ 
($\alpha=1 \sim n$) and their derivatives $\partial_\alpha$ are obtained 
as follows,
$$ 
r\theta^\alpha \theta^\beta=-\theta^\beta \theta^\alpha, \quad 
\partial_\alpha \partial_\beta=-r \partial_\beta \partial_\alpha, 
\qquad (\alpha \leq \beta),$$
$$ r\partial_\alpha \theta^\beta = - \theta^\beta \partial_\alpha, 
\qquad (\alpha \neq \beta), 
\eqno(2.11)$$
$$ \partial_\alpha \theta_\alpha = 1- \theta_\alpha \partial_\alpha 
+(r^{-2}-1)\sum_{\beta > \alpha} \theta ^\beta \partial _\beta.
$$
For all generators $T^\alpha_{\alpha +1}$, $T^{\alpha +1}_\alpha$ and 
$H_\alpha$ of deformed $su(2)$ parts, operations on the above quantum 
space $\theta^\alpha$ and $\partial_\alpha$ are cast into the following 
commutation relations 
$$
[T^\alpha_{\alpha +1},\theta^\alpha]_r=0, \quad 
[T^\alpha_{\alpha +1},\partial_\alpha]_{1/r}=-r^{-1}\partial_{\alpha +1},$$
$$[T^\alpha_{\alpha +1},\theta^{\alpha+1}]_{1/r}=\theta^\alpha, \quad 
[T^\alpha_{\alpha +1},\partial_{\alpha+1}]_r=0,$$
$$[T^{\alpha +1}_\alpha ,\theta^\alpha]_r=\theta^{\alpha+1}, \quad 
[T^{\alpha +1}_\alpha ,\partial_\alpha]_{1/r}=0,$$
$$[T^{\alpha +1}_\alpha ,\theta^{\alpha+1}]_{1/r}=0, \quad 
[T^{\alpha +1}_\alpha ,\partial_{\alpha+1}]_r=-r\partial_\alpha,
\eqno(2.12)$$
$$[H_\alpha ,\theta^\alpha]_{r^2}=-r\theta^\alpha, \quad 
[H_\alpha ,\partial_\alpha]_{1/r^2}=r^{-1}\partial_\alpha ,$$
$$[H_\alpha ,\theta^{\alpha+1}]_{1/r^2}=r^{-1}\theta^{\alpha +1}, \quad 
[H_\alpha ,\partial_{\alpha+1}]_{r^2}=-r\partial_{\alpha +1}.
$$
These generators commute with the other coordinates and derivatives. 

Hence, a deformed $su(n)$ algebra as an abstract algebraic relation is 
obtained as follows: for the $su(2)$ parts, 
$$
[H_\alpha,T^\alpha_{\alpha +1}]_{r^4}=-r^2(r+r^{-1})
T^\alpha_{\alpha +1},\quad
[H_\alpha,T^{\alpha +1}_\alpha ]_{1/r^4}=r^{-2}(r+r^{-1})
T^{\alpha +1}_\alpha ,$$
$$[T^{\alpha +1}_\alpha,T^\alpha_{\alpha +1}]_{r^2}=rH_\alpha,
\eqno(2.13)
$$
and the other generators $T^\alpha_\gamma$ of the deformed $su(n)$ are  
as \cite{Kobayashi}
$$
T^\alpha_\gamma \equiv [T^\alpha_\beta,T^\beta_\gamma]_r, \qquad 
(\alpha < \beta < \gamma {\rm \ or \ } \alpha > \beta > \gamma).
\eqno(2.14)
$$
Using (2.12), we can obtain operation of $T^\alpha_\gamma$ on the quantum 
space. Also the whole algebra of the deformed $su(n)$ can be obtained 
\cite{Kobayashi}.

In closing the section, we show a realization of all the generators 
$T^\alpha_{\alpha +1}$, $T^{\alpha +1}_\alpha$ and $H_\alpha$ which satisfy 
the commutation relations (2.13) and (2.14) using the above $n$-dimensional 
quantum space as below. Differently from two-dimensional case, it is somewhat 
complicated rather than (2.10). Here we define $\Theta^\alpha$, 
$\Theta^{\alpha +1}$, $D_\alpha$  and $D_{\alpha +1}$ as 
$$
\Theta^\alpha \equiv \mu'_\alpha \mu_{\alpha}^{-1/2} \theta^\alpha, 
\quad 
D_\alpha \equiv \mu'_\alpha \mu_{\alpha}^{-1/2} \partial_\alpha,$$
$$\Theta^{\alpha +1} \equiv \mu_{\alpha +1}^{-1/2} \theta^{\alpha +1}, 
\quad 
D_{\alpha +1} \equiv \mu_{\alpha +1}^{-1/2} \partial_{\alpha +1},
\eqno(2.15)
$$
where $\mu'_\alpha$ and $\mu_\alpha$ are defined as 
\cite{Ogievetsky,Kobayashi}
$$ 
\mu'_\alpha \equiv [1+(r^{-2}-1)\theta^{\alpha +1} \partial_{\alpha +1}
]^{1/2}, \quad 
 \mu_\alpha \equiv 1 +(r^{-2}-1)\sum_{\beta > \alpha} \theta^\beta 
\partial_\beta.
\eqno(2.16)
$$
It is easy to show that these coordinates and derivatives, $\Theta^\alpha$, 
$\Theta^{\alpha +1}$, $D_\alpha$ and $D_{\alpha +1}$, satisfy the same 
commutation relations as the two-dimensional quantum space.Using them, we can 
realize the generators, $T^\alpha_{\alpha +1}$, $T^{\alpha +1}_\alpha$ and 
$H_\alpha$ as 
$$
T^\alpha_{\alpha +1}=\Theta^\alpha D_{\alpha +1}, \quad 
T^{\alpha +1}_\alpha=\Theta^{\alpha +1} D_{\alpha} $$
$$H_\alpha=r^{-1}\Theta^{\alpha +1} D_{\alpha +1}-r\Theta^\alpha D_\alpha.
\eqno(2.17)
$$

In the case of quantum bosonic spaces, we can construct the same deformed 
$su(n)$ algebra in the same way as explained above. The operation rule of 
the generators on the bosonic quantum space is same as (2.12). However, we 
can not realize the deformed $su(n)$ algebra in a simple way as (2.10) or 
(2.17) using coordinates and derivatives of the quantum bosonic space.

%%%%%%%%%%%%%%%%%%%%%%%% Section 3 %%%%%%%%%%%%%%%%%%%%%%%%%%%%%%%%%%%%%%%%%%
\vspace{0.8 cm}
\leftline{\large \bf 3. Quantum deformed $gl(n)$ and $igl(n)$ algebras}
\vspace{0.8 cm}
%%%%%%%%%%%%%%%%%%%%%%%%%%%%%%%%%%%%%%%%%%%%%%%%%%%%%%%%%%%%%%%%%%%%%%%%%%%%%%

In this section, we define a dilatation generator $D$ in order to construct 
quantum deformed $gl(n)$ and $igl(n)$ algebras. Let us start with the case 
of the 2-dimensional quantum fermionic space. As a first step, we utilize 
$\theta^1 \partial_1$ and $\theta^2 \partial_2$ in order to realize a 
generator $D$ which decouples to the three generators of (2.10). Namely, we 
define the generator $D$ as 
$$
D \equiv r^{-1}(\theta^1 \partial_1+\theta^2 \partial_2).
\eqno(3.1)
$$
We see that the above generator operates on the quantum space as
$$
[D,\theta^\alpha ]_{1/r^2}=r^{-1}\theta^\alpha, \quad 
[D,\partial_\alpha ]_{r^2}=-r\partial_\alpha, 
\eqno(3.2)
$$
where $\alpha=1,2$ and thus we can show that $D$ commutes with the three 
generators of (2.10), even in the case with the two deformation parameters.
We therefore understand that the deformed $gl(2)$ algebra is obtained in 
terms of $T^1_2$, $T^2_1$, $H_1$ and $D$. Even though we have considered the 
quantum space with two deformation parameters, this deformed $gl(2)$ algebra 
depends only on one parameter $r$. Another deformation parameter appears when 
we consider an inhomogeneous extension. In addition to these generators, 
$\partial_1$, $\partial_2$ define a deformed $igl(2)$ algebra. The second 
deformation parameter $q$ appears in the commutation relations between 
$T^\alpha_\beta$ and $\partial_\alpha$ as well as the commutation relation 
between $\partial_1$ and $\partial_2$.

In general, we can construct the deformed $gl(n)$ and $igl(n)$ algebras 
similarly. Here we define the dilatation generator $D$ as follows,
$$
D \equiv r^{-1}\sum ^n_{\alpha=1} \theta^\alpha \partial_\alpha.
\eqno(3.3)
$$
The generator $D$ commutes with all the generators of the deformed $su(n)$ 
algebra, $T^\alpha_\gamma$ and $H_\alpha$. These generators span the deformed 
$gl(n)$ algebra. Moreover $D$ operates on the quantum space as 
$$
[D,\theta^\alpha ]_{1/r^2}=r^{-1}\theta^\alpha, \quad 
[D,\partial_\alpha ]_{r^2}=-r\partial_\alpha.
\eqno(3.4)
$$
We can construct the deformed $igl(n)$ algebra in terms of $T^\alpha_\gamma$, 
$H_\alpha$, $D$ and $\partial_\alpha$. In other words, all the operators 
$\theta^\alpha \partial_\gamma$ of the $n$-dimensional quantum space 
construct the deformed $gl(n)$ algebra. In addition to these, the derivatives 
$\partial_\alpha$ derive the deformed $igl(n)$ algebra. It is easy to extend 
the construction of $igl(n)$ to cases with multi-parameters.

Finally, we comment on the construction on the quantum bosonic space. In 
this case, we can define $D$ which has the same commutation relations as 
(3.4). We can construct the deformed $igl(n)$ algebra in terms of $T^i_j$, 
$H_i$, $D$ and $\partial / \partial x_i$, although these deformed $gl(n)$ 
generators can not be realized by $x^i$ and 
$\partial / \partial x_i$ as (2.10).

%%%%%%%%%%%%%%%%%%%%%%%%%% Section 4 %%%%%%%%%%%%%%%%%%%%%%%%%%%%%%%%
\vspace{0.8 cm}
\leftline{\large \bf 4. Map to Castellani's deformation}
\vspace{0.8 cm}
%%%%%%%%%%%%%%%%%%%%%%%%%%%%%%%%%%%%%%%%%%%%%%%%%%%%%%%%%%%%%%%%%%%%%

We here present some mapping relations from our deformation (for simplicity, 
we analyze one parameter case substituting $r$ with $q$) to another 
deformation which possesses a non-trivial Hopf algebra structure discussed 
by Castellani \cite{Castellani}. There is a difference between $gl$ and $igl$ 
cases. To show it, we first consider a map to the following algebraic 
relation for deformed $gl(2)$ (hereafter we call it $\hat\chi$-algebra), 
$$
{\hat\chi_+}{\hat\chi_1}-r^{-2}{\hat\chi_1}{\hat\chi_+}
={\hat\chi_+},\quad 
{\hat\chi_1}{\hat\chi_-}-r^{-2}{\hat\chi_-}{\hat\chi_1}={\hat\chi_-},$$
$$ r^2{\hat\chi_+}{\hat\chi_-}-r^{-2}{\hat\chi_-}{\hat\chi_+}
={\hat\chi_2}-{\hat\chi_1},\quad
{\hat\chi_\pm}X-r^{\pm 2}X{\hat\chi_\pm}=0, $$
$$[{\hat\chi_1},{\hat\chi_2}]=0, \quad X=(1-r^2)X{\hat\chi_2}-r.\eqno(4.1)
$$
The deformed $gl(2)$ of Castellani is obtained through the rescaling 
$\chi_i=X{\hat\chi_i}$ ($i=1,2,+,-$). We introduce the following 
recombination for convenience,
$$
A={q\over q+q^{-1}}(qD+q^{-1}H),\quad B={1\over q+q^{-1}}(D-H),\eqno(4.2)
$$
and consider the case of ${\hat\chi_+}=T^2_1$ and ${\hat\chi_-}=T^1_2$. 
Expanding $\chi_1,\chi_2$ and $X$ as linear combinations of our basis 
$A,B$ and $AB$, we find a family of solutions which satisfy the 
$\hat\chi$-algebra (4.1) as follows 
$$ 
   \hat\chi_2 = aA + bB + cAB, $$
$$ \hat\chi_1 = (a-r^2)A + (b+r^{-2})B + (c+r^2q^2-r^{-2})AB, \eqno(4.3a)$$
$$ X = (x+r)A + (y+r)B + (z-r-y)AB - r, $$
where the parameters $a,b,c,x,y,z$ are arbitrary as long as they satisfy 
$$
b={a-1\over r^2},\quad x={r\over 1-(1-r^2)a},\quad 
y={r^3\over 1-(1-r^2)a},$$
$${z-r-y\over 1-r^2}=z(a+q^{-2}(b+c)) + (x-y)b + xc.\eqno(4.3b) $$
If we fix the parameter $a$, the degree of freedom of undetermined 
parameters is two because $b,x,y$ are also fixed in conformity with the 
value of $a$. This means that we have a two-parameter family of relations 
between different deformation parameters $q$ and $r$ of each deformed $gl(2)$. 
Similarly, we can observe another family of mapping for the choice 
$\chi_+\propto T^1_2$ and $\chi_-\propto T^2_1$. Note that the above case is 
$\chi_+=xT^2_1$ and $\chi_-=yT^1_2$. For example, choosing 
$\chi_+=r^2T^1_2$ and $\chi_-=r^2T^2_1$, we obtain
$$
\chi_1=r^3B + \alpha AB,\quad \chi_2= rA + \beta AB,\eqno(4.4a)$$
with
$$
(1-r^2)\left(r\beta + (r+\beta q^{-2})(\beta-\alpha-r^3)\right)
+r(\alpha-\beta)=r^4(q^2-1).\eqno(4.4b)
$$
Although we have fixed some coefficients in the second family (4.4), 
we still have unknown parameters $\alpha$ and $\beta$. Namely, as for 
isomorphisms on deformed $gl(2)$ from our deformation to Castellani's one, 
we have a relation between deformation parameters which is expressed by 
the curves with two parameters as (4.3b) or (4.4b). 

Next, we consider isomorphisms on inhomogeneous extensional part. 
In this case, all the redundant free parameters are fixed and the relation 
between deformation parameters become a straight-line without free 
parameter. Following Castellani's deformation, we have two deformation 
parameters $r$ and $s$ in their algebra,
$$  r^2\chi_1 P_1  -  P_1 \chi_1  =  (1-r^2) P_2\chi_-    -   r P_1, $$
$$  r^2\chi_2 P_2  -  P_2 \chi_2  =                       -   r P_2, $$
$$  s  \chi_+ P_1  -  P_1 \chi_+  =  r^2(r^2-1) \chi_2P_2  -  r^3 P_2, $$
$$  r^2\chi_- P_2  -  sP_2\chi_-  =  s(r^2-1)\chi_2P_1     -  sr P_1,
\eqno(4.5) $$
$$  s  \chi_+ P_2  -  P_2 \chi_+  =  0,\quad\,
    r^2\chi_- P_1  -  s P_1\chi_- =  0, $$
$$  [\,\chi_2\,,P_1\,] =0, \,\quad\, [\,\chi_1\,,P_2\,] = (1-r^2)
{s\over r^2}\chi_+P_1,$$
$$  r^2 P_1P_2  -  s P_1P_2  = 0.$$
Under this extension, there is no solution for arbitrary free parameters 
in both cases of (4.3) and (4.4). We have to fix the values of free 
parameters in order to satisfy (4.5). This seems to uniquely determine the 
values of free parameters, however in addition, we must not forget 
considering a map between our $\partial_i$ and $P_i$. Then we can find a 
solution which satisfies (4.5). We have no solution for the first map 
(4.3), while the second type (4.4) has the solution 
$$ 
             \alpha = \beta = 0,      \eqno(4.6) 
$$
and the equation (4.4b) reduces into
$$ 
                           r = q.   \eqno(4.7) 
$$ 
Mapping from our coordinate derivatives to $P_1$,$P_2$ and another 
deformation parameter $s$ in (4.5) are given by
$$
P_1=\partial_1 + (r^2-1)A\partial_1,\,\quad\, P_2=r^{-1}\partial_2,
\eqno(4.8)
$$
$$             s =q^{-1}.\eqno(4.9) $$
Hence, the equations (4.4) with (4.6)-(4.9) connect our deformed $igl(2)$ 
to (4.1) and (4.5).

As concluding remarks, we notice a duality on $gl(2)$ part of the map 
(4.4) with the case of (4.6). Exchanging $1 \leftrightarrow 2$, 
$+\leftrightarrow-$ and $q\leftrightarrow q^{-1}$, we can transform (4.4) 
into the following another map
$$
\chi_+ = r^2T^2_1, \quad \chi_- = r^2T^1_2, \eqno(4.10a)$$
$$
\chi_1 = r^3A, \quad \chi_2 = rB,\quad \mbox{with}\quad q=1/r.
\eqno(4.10b) $$
As mentioned above, there is no free parameters in this case and so the 
relation between $q$ and $r$ is fixed. It should be noticed that this map 
is not included in (4.3) as a special case.

%%%%%%%%%%%%%%%%%%%%%%%%%% Section 5 %%%%%%%%%%%%%%%%%%%%%%%%%%%%%%%%
\vspace{0.8 cm}
\leftline{\large \bf 5. Conclusion and discussion}
\vspace{0.8 cm}
%%%%%%%%%%%%%%%%%%%%%%%%%%%%%%%%%%%%%%%%%%%%%%%%%%%%%%%%%%%%%%%%%%%%%

In this paper, we have introduced a method how to construct a deformed 
$igl(n)$ algebra based on the fermionic and/or bosonic quantum space. In 
section 2, we have discussed the construction of a deformed $sl(n)$ algebra 
which is independent of realizations. After that, incorporating the dilatation 
and translation operators on the quantum space, we have embodied a deformed 
$igl(n)$ algebra in section 3. In the case of $igl(2)$ deformation, we 
have shown a mapping relation from our deformation to other one 
\cite{Castellani}. 

As a result of the analysis of section 4, there exist at least two 
categories of isomorphisms between our deformed $gl(2)$ to the Castellani 
deformation. The first one, to which (4.3) belongs, does not hold in the case 
of the inhomogeneous extension. The second one, to which (4.4) belongs, 
includes a mapping relation which holds in the inhomogeneous case. As for 
(4.10), we have not examined whether its inhomogeneous case holds or not, 
however it might belong to the second one for the reasons that there exists 
a duality to the (4.4) with (4.6) and that it certainly different from the 
type of (4.3) as previously suggested. It would be interesting to investigate 
complete classification of these kind of maps and a duality relation.

The advantage of our approach consists in simple construction 
of representations of quantum algebras. Owing to this property, 
we have found the relation between our deformation and the other 
which possesses non-trivial Hopf algebra structure. 
It is not so trivial whether or not we can find such map relation in 
other cases. It is also interesting to investigate our approach applying 
to the $iso(n)$, $isp(n)$, Poincar\'e and super algebras.  

\vspace{1cm}
\noindent
{\em Acknowledgment}

We would like to thank T. Uematsu for useful suggestions.
% 
% 
%%%%%%%%%%%%%%%%%%%%%%%%%%%%%%%%%%%%%%%%%%%%%%%%%%%%%%%%%%%%%%%%%%%%%%
%%%%%%%%%%%%%%%%%%%%%%% References %%%%%%%%%%%%%%%%%%%%%%%%%%%%%%%%%%%
%%%%%%%%%%%%%%%%%%%%%%%%%%%%%%%%%%%%%%%%%%%%%%%%%%%%%%%%%%%%%%%%%%%%%%

\end{document}